\documentclass[11pt,reqno]{article}
\usepackage{latexsym, amsmath, amssymb, a4, epsfig}
\usepackage{graphicx}

\newcommand{\nm}{\noalign{\smallskip}}
\newcommand{\ds}{\displaystyle}

\newcommand{\p}{\partial}
\newcommand{\pd}[2]{\frac {\p #1}{\p #2}}
\newcommand{\eqnref}[1]{(\ref {#1})}

\newcommand{\Cbb}{\mathbb{C}}

\newcommand{\Rbb}{\mathbb{R}}

\newcommand{\la}{\langle}
\newcommand{\ra}{\rangle}

\newcommand{\Hcal}{\mathcal{H}}

\newcommand{\Kcal}{\mathcal{K}}
\newcommand{\Ccal}{\mathcal{C}}

\newcommand{\Scal}{\mathcal{S}}

\newcommand{\Ga}{\alpha}
\newcommand{\Gb}{\beta}

\newcommand{\Gvf}{\varphi}
\newcommand{\Gg}{\gamma}

\newcommand{\Gl}{\lambda}

\newcommand{\Gt}{\theta}

\newcommand{\Gs}{\sigma}

\newcommand{\Gz}{\zeta}
\newcommand{\GD}{\Delta}

\newcommand{\GO}{\Omega}

\newcommand{\beq}{\begin{equation}}
\newcommand{\eeq}{\end{equation}}

\numberwithin{equation}{section}
\numberwithin{figure}{section}

\begin{document}
\title{Construction of conformal mappings by generalized polarization tensors\thanks{\footnotesize This work is supported by the Korean Ministry of Education, Sciences and Technology through NRF grants Nos. 2010-0004091, 2010-0017532 and 2013-012931}}

\author{Hyeonbae Kang\thanks{Department of Mathematics, Inha University, Incheon
402-751, Korea (hbkang, hdlee@inha.ac.kr).} \and
Hyundae Lee\footnotemark[2]  \and Mikyoung Lim\thanks{Department of Mathematical Sciences, Korea Advanced Institute of Science and Technology, Daejeon 305-701, Korea (mklim@kaist.ac.kr).}\footnotemark[3]
}

\date{}
\maketitle

\begin{abstract}
We present a new systematic method to construct the conformal mapping from outside the unit disc to outside of a simply connected domain using the generalized polarization tensors. We also present some numerical results to validate effectiveness of the method.
\end{abstract}

%%%%%%%%%%%%%%%%%%%%%%%%%%%%%%%%%%%%%%%%%%%%%%%%%%%%%%%%
\section{Introduction}
%%%%%%%%%%%%%%%%%%%%%%%%%%%%%%%%%%%%%%%%%%%%%%%%%%%%%%%%

Riemann mapping theorem tells us that if the domain $\GO$ is simply connected, then there is a conformal mapping from $\Cbb \setminus U$ ($U$ is the unit disc) onto $\Cbb \setminus \GO$ of the form
\beq\label{conformal}
\Phi(\Gz)= \mu_{-1} \Gz + \mu_0 + \frac{\mu_1}{\Gz} + \frac{\mu_2}{\Gz^2} + \cdots,
\eeq
and the mapping is unique under the assumption $\mu_{-1} >0$. The purpose of this paper is to present a new method to compute the coefficients $\mu_{-1}, \mu_0, \mu_1, \ldots$ of the mapping.

Since the conformal mapping plays a fundamental role in various areas of mathematics and applications, many methods to construct conformal mappings have been introduced, for which we refer readers to \cite{hen} and comprehensive references therein instead of citing a long list of literature on numerical computation of the conformal mapping. The method of this paper uses the generalized polarization tensors (GPTs). The GPT is a sequence of tensors (matrices in two dimensions) associated with a domain which appears naturally in the multi-polar expansion of the electric potential. It contains rich information of the shape of the domain. For example, it is proved in \cite{ak02} that the full set of GPTs determines the domain uniquely. The notion of GPTs has been used in various areas of applications such as inverse problems and imaging and the theory of composites. We refer to \cite{book2, ak11, kang, milton} and references therein for these applications. More recent applications of GPT include shape representations \cite{agkly, aklz}, dictionary matching \cite{abgj, ackw}, invisibility cloaking \cite{akll}, and electro-sensing \cite{abg, abgw}.

In this paper we derive canonical relations between GPTs and coefficients of the conformal mapping. Since GPTs of a domain can be computed numerically using the boundary integral method (see section \ref{sec2}), so can the coefficients of the conformal mapping using these relations. We will show some numerical examples of the ranges of mappings $\mu_{-1} \Gz + \mu_0 + \frac{\mu_1}{\Gz} + \cdots + \frac{\mu_n}{\Gz^n}$ for $n=1, 2, \ldots$. They clearly exhibit how the ranges gradually approximate the given domain.

This paper is organized as follows. In section \ref{sec2} we review the definition and computation of GPT, and its relation to eigenvalues of Neumann-Poincar\'e operator. Section \ref{sec3} is to derive the relation between GPTs and coefficients of the conformal mapping. Some numerical examples are provided in section \ref{sec4}. The paper is concluded with some discussions.

%%%%%%%%%%%%%%%%%%%%%%%%%%%%%%%%%%%%%%%%%%%%%%%%%%%%
\section{GPTs and eigenvalues of Neumann-Poincar\'e operator}\label{sec2}
%%%%%%%%%%%%%%%%%%%%%%%%%%%%%%%%%%%%%%%%%%%%%%%%%%%

Let $\GO$ be a domain with the Lipschitz boundary in $\Rbb^2$ and suppose that the conductivity (or the dielectric constant) of $\GO$ is $k$ and that of the background is $1$ ($k \neq 1$). So, the distribution of the conductivity is given by
\beq
\Gs=k \chi(\GO) + \chi(\Rbb^2 \setminus \overline{\GO}),
 \eeq
where $\chi$ denotes the indicator function. For a given harmonic function $h$ in $\Rbb^2$ we consider the following transmission problem:
\beq\label{trans}
\left\{
\begin{array}{ll}
\nabla \cdot \Gs \nabla u = 0 \quad &\mbox{in } \Rbb^2, \\
\nm
u(x)-h(x) = O(|x|^{-1}) \quad &\mbox{as } |x| \to \infty.
\end{array}
\right.
\eeq

If $h$ takes the form in polar coordinates
\beq
h(x)=a_0+\sum_{n=1}^\infty r^n (a_n^c \cos n \theta + a_n^s \sin n\theta),
\eeq
then it is known \cite{book2} that the solution $u$ to \eqnref{trans} can be represented as
\begin{align}
(u-h)(x) & = -\sum_{m=1}^\infty\frac{\cos m\theta}{2\pi mr^m}\sum_{n=1}^\infty
\bigr(M_{mn}^{cc}a_n^c + M_{mn}^{cs}a_n^s \bigr) \nonumber \\
& \quad -\sum_{m=1}^\infty\frac{\sin m\theta}{2\pi m
r^m}\sum_{n=1}^\infty
\bigr(M_{mn}^{sc}a_n^c + M_{mn}^{ss}a_n^s \bigr) \quad\mbox{as
} |x| \to \infty, \label{expan2}
\end{align}
The quantities $M_{mn}^{\Ga\Gb}$ ($\Ga, \Gb=c,s$) appearing in the expansion \eqnref{expan2} are called (contracted) generalized polarization tensors (GPTs).

We emphasize that GPTs can be computed numerically once the domain is given. In fact, let
\beq
P_n^c(x) = r^n \cos n\Gt \quad \mbox{and} \quad P_n^s(x) = r^n \sin n\Gt.
\eeq
Then $M_{mn}^{\Ga\Gb}$, $\Ga, \Gb=c,s$, are given by
\beq\label{gpt}
M^{\Ga\Gb}_{mn} = \int_{\p \GO} P_m^\Gb(x)  (\Gl
I - \Kcal^*_{\p\GO})^{-1}[\nu \cdot \nabla P_n^\Ga ](x) \, d\Gs,
\eeq
where
\beq
\Gl= \frac{k+1}{2(k-1)},
\eeq
and $\Kcal^*_{\p\GO}$ is the Neumann-Poincar\'e (NP) operator defined by
\beq \label{introkd2}
\mathcal{K}^*_{\p\GO} [\Gvf] (x) = \frac{1}{2\pi} \int_{\p\GO}
\frac{\la x -y, \nu_x \ra}{|x-y|^2} \Gvf(y)\,d\Gs(y)\;, \quad x \in \p\GO.
\eeq
Here $\nu_x$ is the outward unit normal vector to $\p B$ at $x$. See \cite{book2, kang} for derivation of \eqnref{gpt}. We emphasize that $|\Gl| \ge 1/2$.

Let us look into the connection between GPTs and eigenvalues of the NP-operator (the reciprocal of the eigenvalues of the NP-operator are called the Fredholm eigenvalues). The connection between Fredholm eigenvalues and conformal mapping was investigated in \cite{schiffer1, schiffer2}. Let $\Scal_{\p\GO} [\Gvf]$ be the single layer potential of a density function $\Gvf\in L^2(\p\GO)$, namely,
\beq
\Scal_{\p\GO} [\Gvf] (x) := \frac{1}{2\pi} \int_{\p \GO} \ln |x-y| \Gvf (y) \, d\sigma(y) \;, \quad x \in \Rbb^2 .
\eeq
The relation between the boundary value of the single layer potential and the NP-operator is given by the following jump formula:
\beq\label{singlejump}
\frac{\p}{\p\nu} \Scal_{\p\GO} [\Gvf] \big |_- (x) = \biggl( - \frac{1}{2} I + \Kcal_{\p\GO}^* \biggr) [\Gvf] (x),
\quad x \in \p\GO\;.
\eeq
Here, $\frac{\p}{\p\nu}$ denotes the normal derivative and the subscript $-$ indicates the limit from the inside $\GO$.

It is known (see, for example, \cite{kang}) that $-\la \Gvf, \Scal_{\p\GO}[\Gvf] \ra$ is an inner product on $L^2_0(\p\GO)$ which is the space of square integrable functions with the mean zero. Let $\Hcal$ be the Hilbert space $L^2_0(\p\GO)$ equipped with this inner product, and define
\beq
\la \Gvf,\psi \ra_\Hcal:= -\la \Gvf, \Scal_{\p\GO} [\psi] \ra, \quad \Gvf, \psi \in \Hcal.
\eeq
Because of Plemelj's symmetrization
principle (also known as Calder\'on's identity)
\beq
\Scal_{\p\GO} \Kcal^*_{\p\GO} = \Kcal_{\p\GO} \Scal_{\p\GO},
\eeq
the operator $\Kcal^*_{\p\GO}$ is self-adjoint on $\Hcal$.

If $\p\GO$ is $\Ccal^{1, \Ga}$ for some $\Ga >0$, then $\Kcal^*_{\p\GO}$ is compact on $\Hcal$. So, $\Kcal^*_{\p\GO}$ has eigenvalues accumulating to $0$.
Let $\Gl_1, \Gl_2, \ldots$ ($|\Gl_1| \ge |\Gl_2| \ge \ldots$) be eigenvalues of $\Kcal^*_{\p\GO}$ on $\Hcal$ counting multiplicities, and $\Gvf_1, \Gvf_2, \ldots$ be the corresponding (normalized) eigenfunctions. Then $|\Gl_n| < 1/2$ for all $n$ and $\Kcal^*_{\p\GO}$ admits the spectral resolution
\beq\label{specresol1}
\Kcal^*_{\p\GO} = \sum_{j=1}^\infty \Gl_j \Gvf_j \otimes \Gvf_j
\eeq
in $\Hcal$. We emphasize that $\{ \Gvf_j \}$ is a basis for $\Hcal$. Using \eqnref{specresol1}, one can easily obtain that
\beq
M_{mn}^{\Ga\Gb} = \sum_{j=1}^\infty \frac{\la \nu \cdot \nabla P_n^\Ga, \Gvf_j \ra_\Hcal \la P_m^\Gb, \Gvf_j \ra}{\Gl - \Gl_j}.
\eeq
In above the second inner product is the usual inner product on $L^2(\p\GO)$. But since $\Kcal^*_{\p\GO}[\Gvf_j] = \Gl_j \Gvf_j$, we have from \eqnref{singlejump} that
$$
\frac{\p}{\p\nu} \Scal_{\p\GO} [\Gvf_j] \big |_- = \biggl( - \frac{1}{2} I + \Kcal_{\p\GO}^* \biggr) [\Gvf_j] = \biggl( \Gl_j - \frac{1}{2} \biggr) \Gvf_j,
$$
and hence
$$
\Gvf_j = \frac{1}{\Gl_j - \frac{1}{2}} \frac{\p}{\p\nu} \Scal_{\p\GO} [\Gvf_j] \big |_-.
$$
Therefore, we have
$$
\la P_m^\Gb, \Gvf_j \ra = \frac{1}{\Gl_j - \frac{1}{2}} \la P_m^\Gb, \nu \cdot \nabla \Scal_{\p\GO} [\Gvf_j]|_{-} \ra
= \frac{1}{\Gl_j-\frac{1}{2}} \la \nu \cdot \nabla P_m^\Gb, \Gvf_j \ra_\Hcal,
$$
where the last equality follows from the divergence theorem. So we have the following relation between GPTs and eigenvalues of NP-operator:
\beq
M_{mn}^{\Ga\Gb} = \sum_{j=1}^\infty \frac{\la \nu \cdot \nabla P_n^\Ga, \Gvf_j \ra_\Hcal \la \nu \cdot \nabla P_m^\Gb, \Gvf_j \ra_\Hcal}{(\Gl-\Gl_j)(\Gl_j-\frac{1}{2})}.
\eeq
We mention that if $k=0$, then $\Gl= -1/2$.

%%%%%%%%%%%%%%%%%%%%%%%%%%%%%%%%%%%%%%%%%%%%%%%%%%%%%%%%%%%%%%%%%%%
\section{GPTs and conformal mappings}\label{sec3}
%%%%%%%%%%%%%%%%%%%%%%%%%%%%%%%%%%%%%%%%%%%%%%%%%%%%%%%%%%%%%%%%%%%

Suppose now that the inclusion is insulated so that $k=0$. Then, the equation \eqnref{trans} is replaced by
\beq\label{condeqn}
\begin{cases}
\GD u =0 \quad&\mbox{in } \Rbb^2\setminus \overline{\GO},\\
\ds \pd{u}{\nu}=0 &\mbox{on } \p \GO,\\
u(x)- h(x) =O(|x|^{-1}) &\mbox{as }  |x|\to \infty.
\end{cases}
\eeq
Let $u$ be the solution to this equation, and let $H$ be an entire function such that $\Re H=h$, and $U$ be an analytic function in $\Cbb \setminus \overline{\GO}$ such that $\Re U=u$. Then $H$ takes the form
$$
H(z) = \Ga_0 +\sum_{n=1}^\infty \alpha_n z^n,\quad\alpha_n= a_n^c - i a_n^s,
$$
and $U$ takes the form
\beq
U(z)=H(z)-\sum_{m=1}^\infty\frac{\beta_m}{z^m},
\eeq
where
\beq
\beta_m:=\frac{1}{2\pi m } \sum_{n=1}^\infty \left[ ( M_{mn}^{cc} a_n^c + M_{mn}^{cs} a_n^s)+ i (M_{mn}^{sc} a_n^c + M_{mn}^{ss} a_n^s)\right].
\eeq
It is more convenient to write $\beta_m$ as
\beq
\beta_m
= \sum_{n=1}^\infty ( \Gg_{mn}^1 \alpha_n + \Gg_{mn}^2 \overline{\alpha_n})
\eeq
with
\begin{align}
\Gg_{mn}^1:=\frac{1}{4\pi m} (M_{mn}^{cc}-M_{mn}^{ss}+i(M_{mn}^{cs}+M_{mn}^{sc})),\\ \Gg_{mn}^2:=\frac{1}{4\pi m} (M_{mn}^{cc}+M_{mn}^{ss}-i(M_{mn}^{cs}-M_{mn}^{sc})).
\end{align}
Then $U$ can be written as
\beq
U (z) = \Ga_0 +\sum_{n=1}^\infty \left(\alpha_n z^n - \sum_{m=1}^\infty \frac{ \Gg_{mn}^1 \alpha_n + \Gg_{mn}^2 \overline{\alpha_n}}{z^m}  \right).
\eeq
So, we have
\beq\label{U(z)}
U (z) = \Ga_0 +\sum_{n=1}^\infty \left[ a_n^c \left( z^n - \sum_{m=1}^\infty \frac{ \Gg_{mn}^1  + \Gg_{mn}^2 }{z^m}  \right) - i a_n^s \left( z^n - \sum_{m=1}^\infty \frac{ \Gg_{mn}^1  - \Gg_{mn}^2 }{z^m}  \right) \right].
\eeq

One can easily see from the Cauchy-Riemann equation that the boundary condition $\pd{u}{\nu}=0$ in \eqnref{condeqn} is equivalent to
\beq
\Im U=\mbox{constant} \quad \mbox{on } \p\GO.
\eeq
Since this condition holds for any entire function $H$, we infer from \eqnref{U(z)} that
\beq\label{constcond}
\Re \left( z^n - \sum_{m=1}^\infty \frac{ \Gg_{mn}^1 -\Gg_{mn}^2 }{z^m}  \right)=\mbox{const.},\quad
\Im \left( z^n - \sum_{m=1}^\infty \frac{ \Gg_{mn}^1 +\Gg_{mn}^2 }{z^m}  \right)=\mbox{const.}
\eeq
on $\p\GO$ for every positive integer $n$.

Let $z=\Phi(\Gz)$ be the conformal mapping from $|\Gz|>1$ onto $\mathbb{C}\setminus \overline{\GO}$, given by \eqnref{conformal}. Let us write $c=\mu_{-1}$ for ease of notation.
Then $U \circ \Phi(\Gz)$ is analytic in $|\Gz|>1$ and takes the form
\begin{align}
U \circ \Phi (\Gz)
&=\Ga_0 +\sum_{n=1}^\infty \left(\alpha_n \Phi(\Gz)^n - \sum_{m=1}^\infty \frac{ \Gg_{mn}^1 \alpha_n + \Gg_{mn}^2 \overline{\alpha_n}}{\Phi(\Gz)^m}  \right).
\end{align}
We infer from \eqnref{constcond} that
\beq \label{constant02}
\Re \left( \Phi(\Gz)^n - \sum_{m=1}^\infty \frac{ \Gg_{mn}^1 -\Gg_{mn}^2 }{\Phi(\Gz)^m}  \right) = \mbox{constant on } |\Gz|=1
\eeq
and
\beq \label{constant01}
\Im \left( \Phi(\Gz)^n - \sum_{m=1}^\infty \frac{ \Gg_{mn}^1 +\Gg_{mn}^2 }{\Phi(\Gz)^m}  \right) = \mbox{constant on } |\Gz|=1.
\eeq
These conditions implies that
\beq \label{f_conjg}
\Phi(\Gz)^n - \sum_{m=1}^\infty \frac{ \Gg_{mn}^1}{\Phi(\Gz)^m} +\overline{\sum_{m=1}^\infty \frac{ \Gg_{mn}^2 }{\Phi(\Gz)^m}}
\eeq
is constant on $|\zeta|=1$.  Since $\sum_{m=1}^\infty \frac{ \Gg_{mn}^2 }{\Phi(\Gz)^m}$ can be expanded as
\beq
\sum_{m=1}^\infty \frac{ \Gg_{mn}^2 }{\Phi(\Gz)^m}= \sum_{k=1}^\infty \frac{s_n}{\zeta^n}, \quad |\zeta|>1,
\eeq
it follows from \eqnref{f_conjg} that
\beq
\Phi(\Gz)^n - \sum_{m=1}^\infty \frac{ \Gg_{mn}^1}{\Phi(\Gz)^m} = \mbox{constant} -  \sum_{k=1}^\infty \bar s_n \zeta^n, \quad |\zeta|=1,
\eeq
so  that
\beq \label{entire01}
\Phi(\Gz)^n - \sum_{m=1}^\infty \frac{ \Gg_{mn}^1  }{\Phi(\Gz)^m} \quad\mbox{is an entire function}.
\eeq

We now derive relations among GPTs and coefficients of the conformal mapping using \eqnref{constant01} and \eqnref{entire01}. We observe, for $\Gz$ with large modulus,
\begin{align}
\frac{1}{\Phi(\Gz)} & = \frac{1}{c\Gz}\cdot\frac{1}{1 + \frac{\mu_0}{c\Gz} + \frac{\mu_1}{c\Gz^2} + \frac{\mu_2}{c\Gz^3} + \ldots} \nonumber\\
& =\frac{1}{c\Gz}\sum_{j=0}^\infty (-1)^j\left(\frac{\mu_0}{c\Gz} + \frac{\mu_1}{c\Gz^2} + \frac{\mu_2}{c\Gz^3} + \ldots\right)^j
= \sum_{k=1}^\infty \frac{B_k}{\Gz^k}, \label{Bk}
\end{align}
where $B_1=1/c$ and
\beq\label{B_k}
B_k:=\frac{1}{c}\sum_{s_1k_1+\ldots+s_jk_j = k-1 \atop s_1,\ldots,s_j > 0,~  k_j> \ldots >k_1 > 0} \left(\frac{-1}{c}\right) ^{s_1+\ldots+s_j}\frac{(s_1+\ldots+s_j)!}{s_1!\cdots s_j!} {\mu_{k_1-1}^{s_1}\cdots\mu_{k_j-1}^{s_j}},\quad k\ge 2.
\eeq
We emphasize that $B_k$ ($k \ge 2$) is determined by $\mu_\ell$ for $\ell \le k-2$. It is helpful to write down first few terms:
\begin{align}
B_2&=-\mu_0/c^2,\nonumber\\
B_3&= -\mu_1/c^2 + \mu_0^2/c^3,\nonumber\\
B_4&= -\mu_2/c^2 +2\mu_0\mu_1/c^3-\mu_0^3/c^4, \label{Bkfewterms} \\
B_5&=-\mu_3/c^2+\mu_1^2/c^3+2\mu_0\mu_2/c^3- 3 \mu_0^2\mu_1/c^4+\mu_0^4/c^5, \nonumber \\
B_6&= -\mu_4/c^2+ 2\mu_0 \mu_3/c^3 + 2\mu_1\mu_2/c^3- 3 \mu_0\mu_1^2/c^4 - \mu_0^5/c^6. \nonumber
\end{align}

We now consider the conditions \eqnref{entire01} when $n=1$.  One can see from \eqnref{Bk} that
$$
\sum_{m=1}^\infty \frac{ \Gg_{m1}^1 }{\Phi(\Gz)^m}
= \sum_{\ell=1}^\infty \sum_{s_1n_1+\ldots+s_jn_j =\ell \atop s_1,\ldots,s_j > 0,~  n_j> \ldots>n_1 > 0} \Gg_{s_1+\ldots+s_j,1}^1 \frac{(s_1+\ldots+s_j)!}{s_1!\cdots s_j!} \frac{B_{n_1}^{s_1}\cdots B_{n_j}^{s_j}}{ \Gz^\ell},
$$
and hence
\begin{align}
&\Phi(\Gz) - \sum_{m=1}^\infty \frac{ \Gg_{m1}^1 }{\Phi(\Gz)^m} \nonumber \\
&=c\Gz + \sum_{\ell=0}^\infty \frac{\mu_\ell}{\Gz^\ell}- \sum_{\ell=1}^\infty \frac{1}{ \Gz^\ell} \sum_{s_1n_1+\ldots+s_jn_j =\ell \atop s_1,\ldots,s_j > 0,~  n_j> \ldots>n_1 > 0} \Gg_{s_1+\ldots+s_j,1}^1 \frac{(s_1+\ldots+s_j)!}{s_1!\cdots s_j!} B_{n_1}^{s_1}\cdots B_{n_j}^{s_j}.
\end{align}
It then follows from \eqnref{entire01} that
\beq \label{mu_higher}
\mu_\ell = \sum_{s_1n_1+\ldots+s_jn_j =\ell \atop s_1,\ldots,s_j > 0,~  n_j> \ldots>n_1 > 0} \Gg_{s_1+\ldots+s_j,1}^1 \frac{(s_1+\ldots+s_j)!}{s_1!\cdots s_j!} B_{n_1}^{s_1}\cdots B_{n_j}^{s_j},\quad \ell \ge 1.
\eeq
We note that $\mu_\ell$ is determined by $\Gg_{m1}^1$ for $m \le \ell$ and $B_k$ for $k \le \ell$ which in turn determined by $\mu_j$ for $j \le k-2$ as we have seen it in \eqnref{B_k}. So, $\mu_\ell$ ($\ell \ge 1$) is determined by $\Gg_{m1}^1$ for $m \le \ell$ and $\mu_0, \ldots, \mu_{\ell-2}$ ($\mu_{-1}=c$).
For example, we have first few terms as follows:
\begin{align}
\mu_1 &= \Gg_{11}^1 B_1, \nonumber \\
\mu_2 &= \Gg_{21}^1 B_1^2 + \Gg_{11}^1 B_2, \nonumber \\
\mu_3 &= \Gg_{31}^1 B_1^3 + 2 \Gg_{21}^1 B_1 B_2 + \Gg_{11}^1 B_3, \label{mufew} \\
\mu_4 &= \Gg_{41}^1 B_1^4 + \Gg_{21} B_2^2 + 3 \Gg_{31}^1 B_1^2 B_2 + 2 \Gg_{21}^1 B_1 B_3 + \Gg_{11}^1 B_4 , \nonumber \\
\mu_5 &= \Gg_{51}^1 B_1^5 + 3\Gg_{31}^1 B_1 B_2^2 + 2 \Gg_{21} (B_1 B_4 + B_2 B_3) + \Gg_{11}^1 B_5, \nonumber \\
\mu_6 &= \Gg_{61}^1 B_1^6 + \Gg_{51}^1 B_1^4 B_2 + \Gg_{41}^1 B_1^3 B_3 + \Gg_{31}^1 (B_1^4 B_4 + B_1 B_2 B_3) + \Gg_{21}^1 (B_1 B_5 + B_2 B_4). \nonumber
\end{align}

We now look into the condition \eqnref{constant01} for $n=1$. One can check that
\begin{align}
&\Phi(\Gz) - \sum_{m=1}^\infty \frac{ \Gg_{m1}^1-\Gg_{m1}^2 }{\Phi(\Gz)^m}\nonumber\\
&=c\Gz + \sum_{\ell=0}^\infty \frac{\mu_\ell}{\Gz^\ell}- \sum_{\ell=1}^\infty \sum_{s_1n_1+\ldots+s_jn_j =\ell \atop s_1,\ldots,s_j > 0,~  n_j> \ldots >n_1 > 0} (\Gg_{s_1+\ldots+s_j,1}^1-\Gg_{s_1+\ldots+s_j,1}^2) \frac{(s_1+\ldots+s_j)!}{s_1!\cdots s_j!} \frac{B_{n_1}^{s_1}\cdots B_{n_j}^{s_j}}{ \Gz^\ell}\nonumber\\
& = c \Gz + \mu_0 +\frac{\mu_1 -(\Gg_{11}^1 -\Gg_{11}^2)/c}{\Gz}
+\frac{\mu_2 +\mu_0(\Gg_{11}^1 -\Gg_{11}^2)/c^2 -(\Gg_{21}^1 -\Gg_{21}^2)/c^2}{\Gz^2}+\ldots
\end{align}
Then using \eqnref{constant01} we obtain
\beq \label{c_mu0}
c^2=-\Gg_{11}^2,\quad  \mu_0=-c^{-2}\Gg_{21}^2,
\eeq
and
\beq \label{vanishing_eq}
\sum_{s_1n_1+\ldots+s_jn_j =\ell \atop s_1,\ldots,s_j > 0,~  n_j>,\ldots>n_1 > 0} \Gg_{s_1+\ldots+s_j,1}^2 \frac{(s_1+\ldots+s_j)!}{s_1!\cdots s_j!} B_{n_1}^{s_1}\cdots B_{n_j}^{s_j}=0,\quad \ell \ge 2.
\eeq

So we conclude that all the coefficients of the conformal mapping is determined from
\beq\label{full_set}
\Gg_{11}^2, ~\Gg_{21}^2,~\{\Gg_{m1}^1\}_{m\in\mathbb{N}}.
\eeq
In fact, $\mu_\ell$ can be determined inductively using these GPTs: $\mu_{-1}=c$ and $\mu_0$ are determined by the formula \eqnref{c_mu0}, $\mu_1$ is determined by the first equation in \eqnref{mufew}, $\mu_\ell$ for $\ell \ge 2$ is determined by formula \eqnref{B_k} and \eqnref{mu_higher} in terms of $\Gg_{m1}^1$ for $m \le \ell$ and $\mu_k$ for $k \le \ell-2$.

%%%%%%%%%%%%%%%%%%%%%%%%%%%%%%%%%%%%%%%%%%%%%%%%%%%%%
\section{Numerical illustration}\label{sec4}
%%%%%%%%%%%%%%%%%%%%%%%%%%%%%%%%%%%%%%%%%%%%%%%%%%%%%

In this section we provide numerical examples of conformal mapping \eqnref{conformal} to outside of simply connected domains obtained using the method presented in the previous section. In order to acquire the GPTs, we solve the boundary integral equation \eqnref{gpt} numerically. We refer readers to \cite{book_new} for more details of the computation and numerical codes.  The number of nodal points used on $\partial\Omega$ is 3072 in each example.

Once GPTs of the given domain are computed, then the first two coefficients $\mu_{-1}$ and $\mu_0$ of the conformal mapping $\Phi$ are determined by \eqnref{c_mu0}, and those of higher order terms by \eqnref{mu_higher}.  Let $\Phi_N$, $N\geq 1$, be the truncation of $\Phi$ at the $N$-th order, namely,
 \beq\label{conformal_M}
\Phi_N(\Gz)= \mu_{-1} \Gz + \mu_0 + \frac{\mu_1}{\Gz} + \frac{\mu_2}{\Gz^2} + \cdots+\frac{\mu_{N}}{\Gz^{N}}.
\eeq

In the following examples, we show the images (in black curves) of the unit circle $(S^1)$ under the $\Phi_N$ for domains $\Omega$ of various shapes. The gray curves are actual boundaries of the domains.

\begin{figure}[ht!]
\begin{center}
\epsfig{figure = 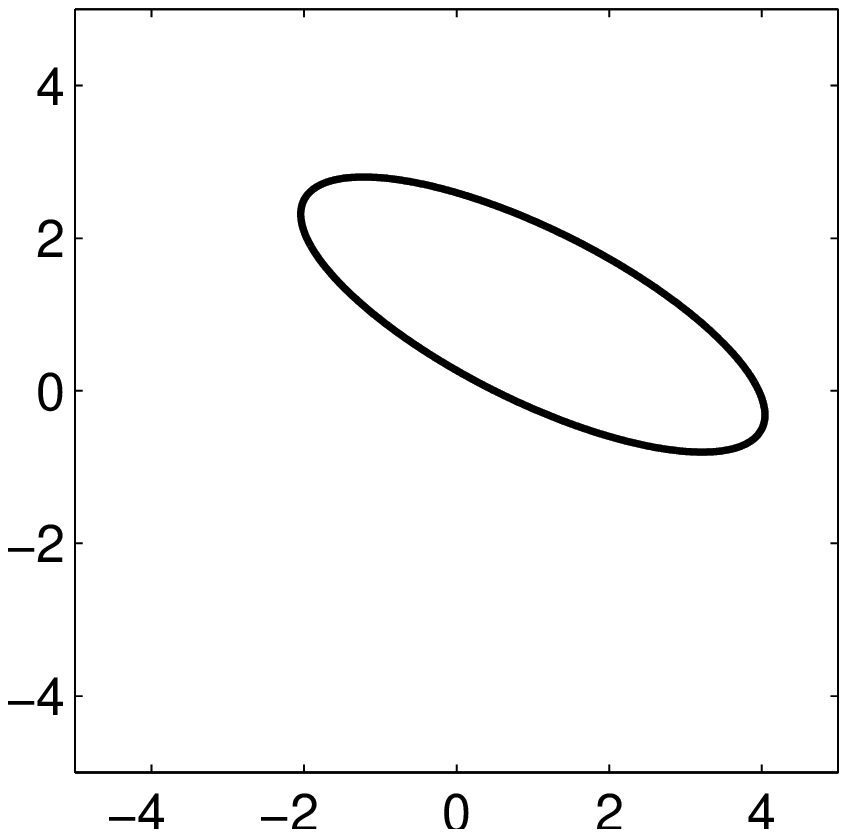,width=3.5cm}\hskip .5cm
\epsfig{figure = 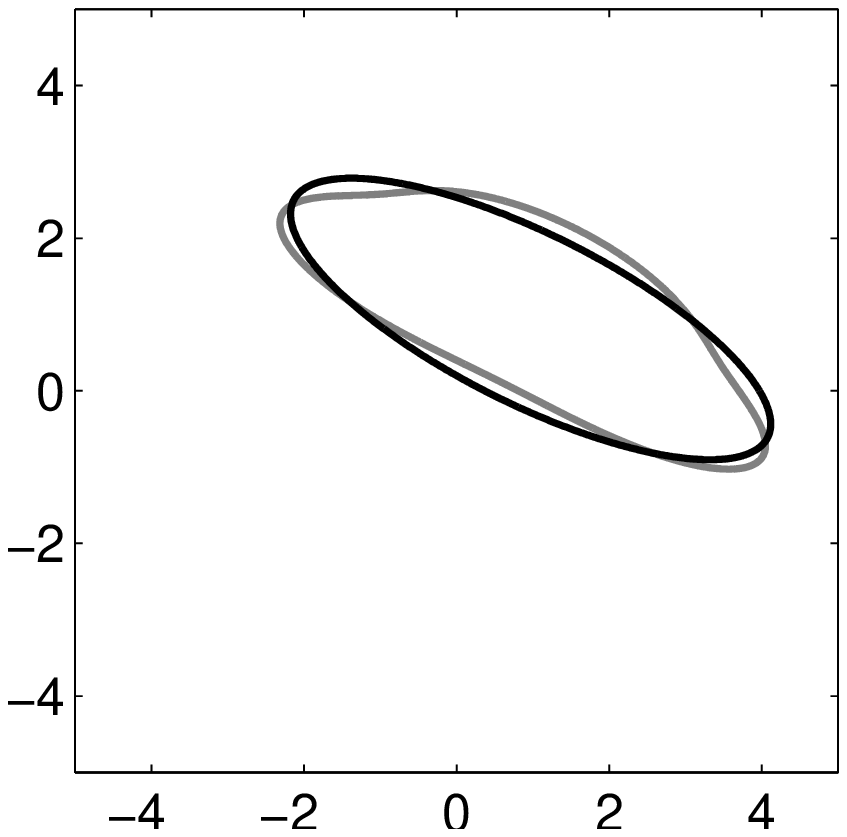,width=3.5cm}\hskip .5cm
\epsfig{figure = 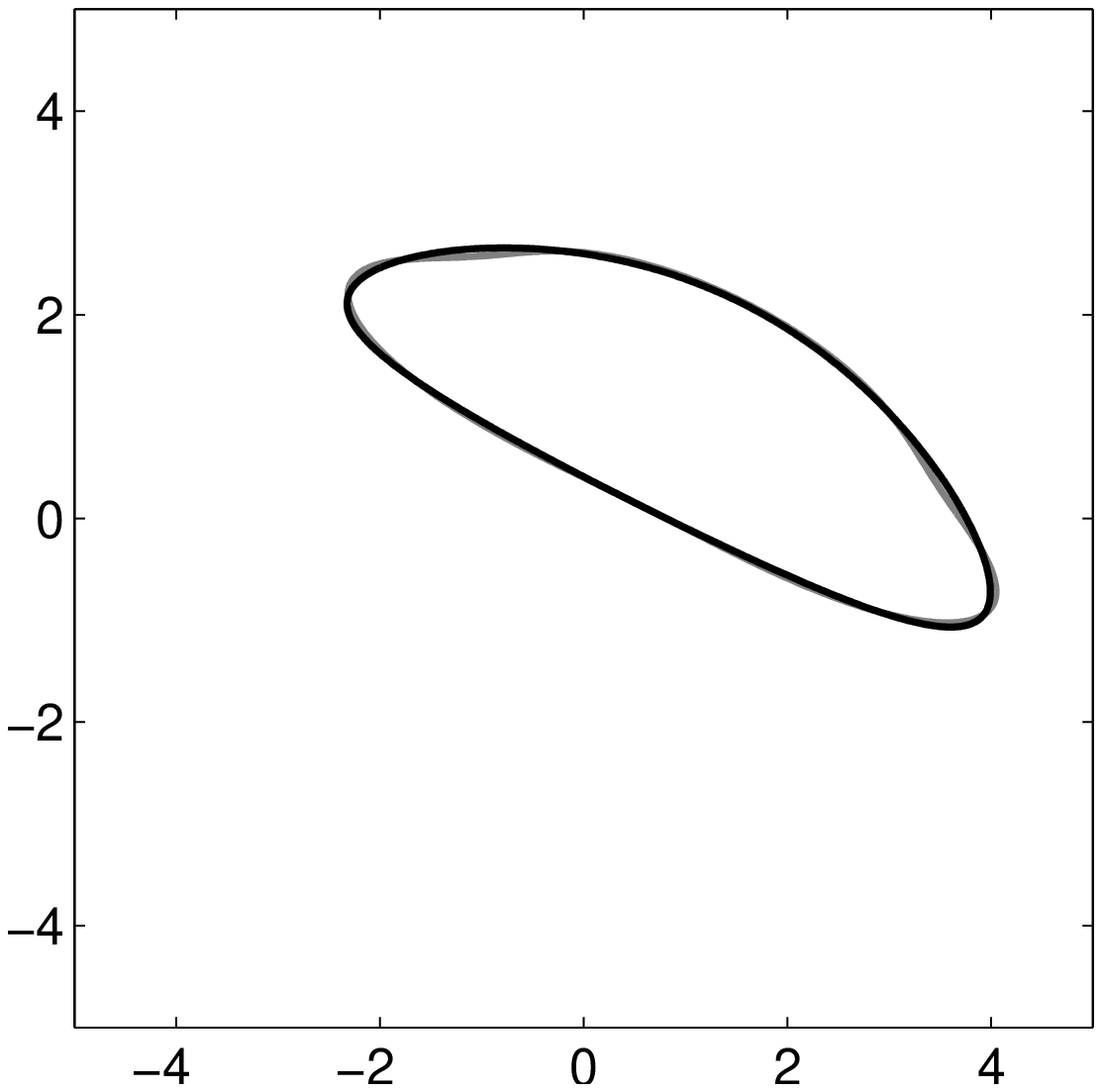,width=3.5cm}
\end{center}
\caption{In the first figure, $\Omega$ is an ellipse and $N=1$. In the next two figures, $\Omega$ is a perturbed ellipse,  and $N$ is 1 and 2 in the middle and the right figures, respectively.}\label{diskellipse}
\end{figure}

\noindent{\bf Example 1}. For ellipses, $\Phi_1(S^1)$ exactly matches with the boundary of $\Omega$. See Figure \ref{diskellipse}. For a perturbed ellipse, $\Phi_N$ with $N=2$ recovers a good approximation of $\p\GO$.

\begin{figure}[ht!]
\begin{center}
\epsfig{figure = 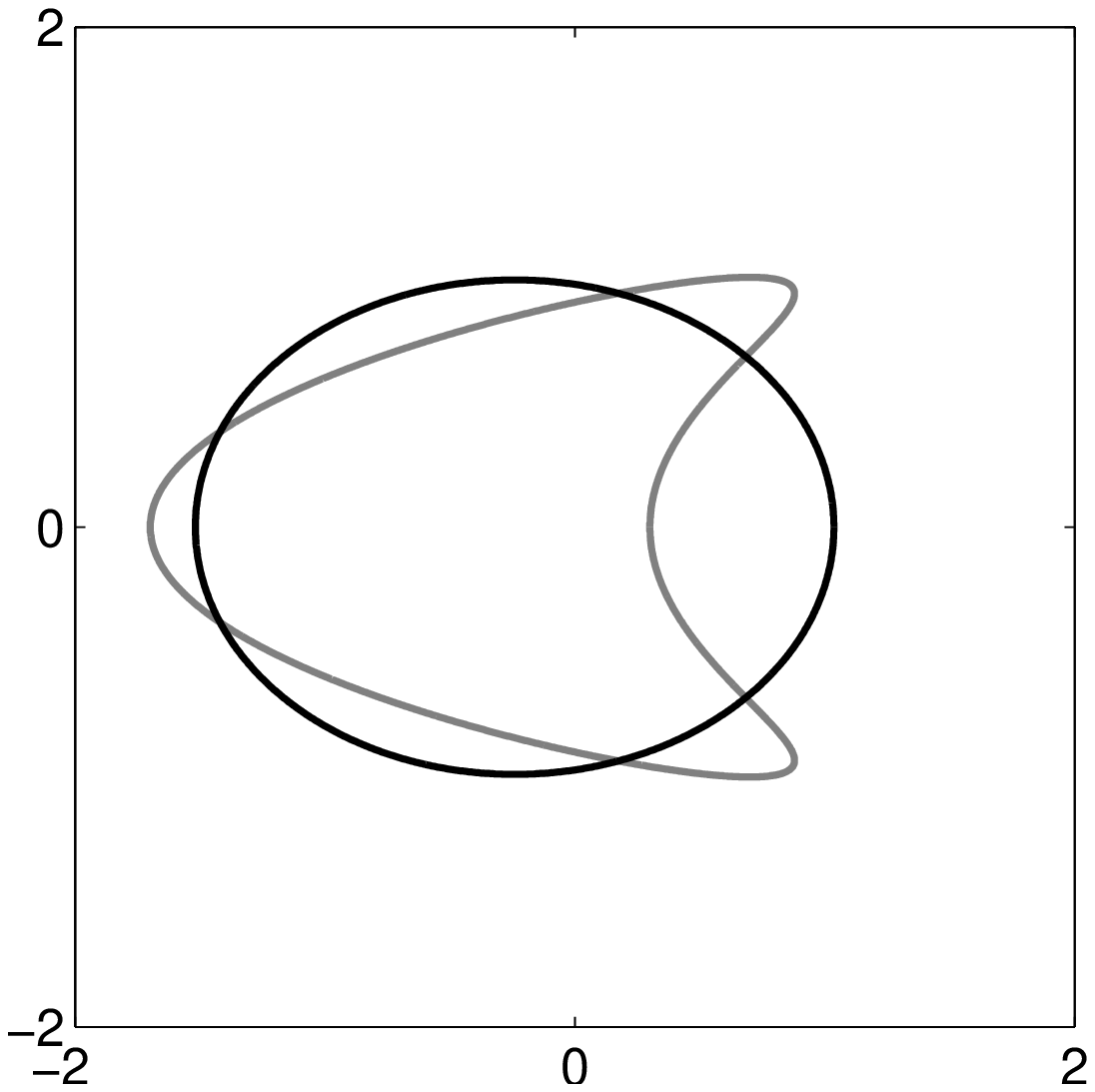,width=3.5cm}\hskip .5cm
\epsfig{figure = 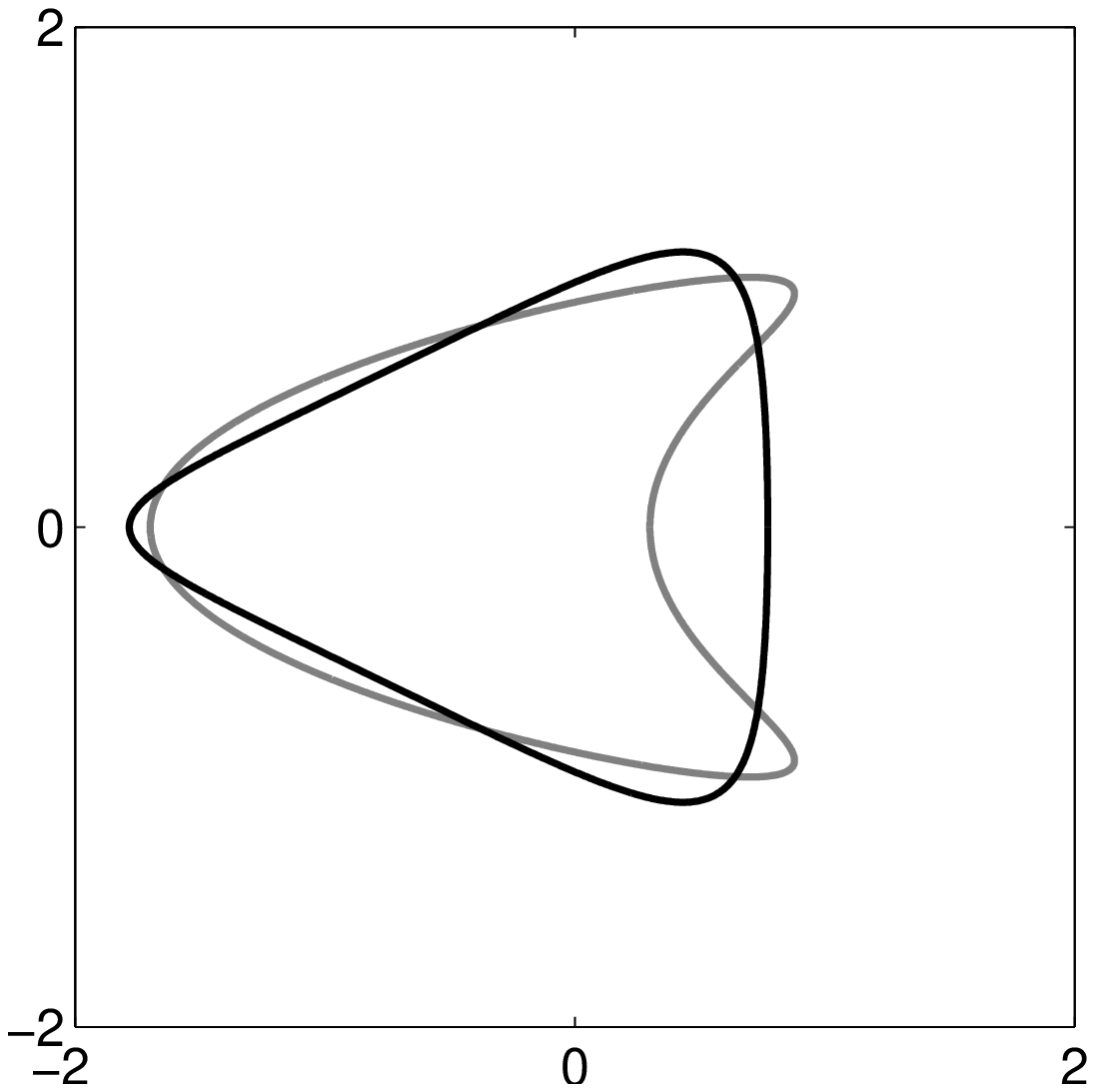,width=3.5cm}\hskip .5cm
\epsfig{figure = 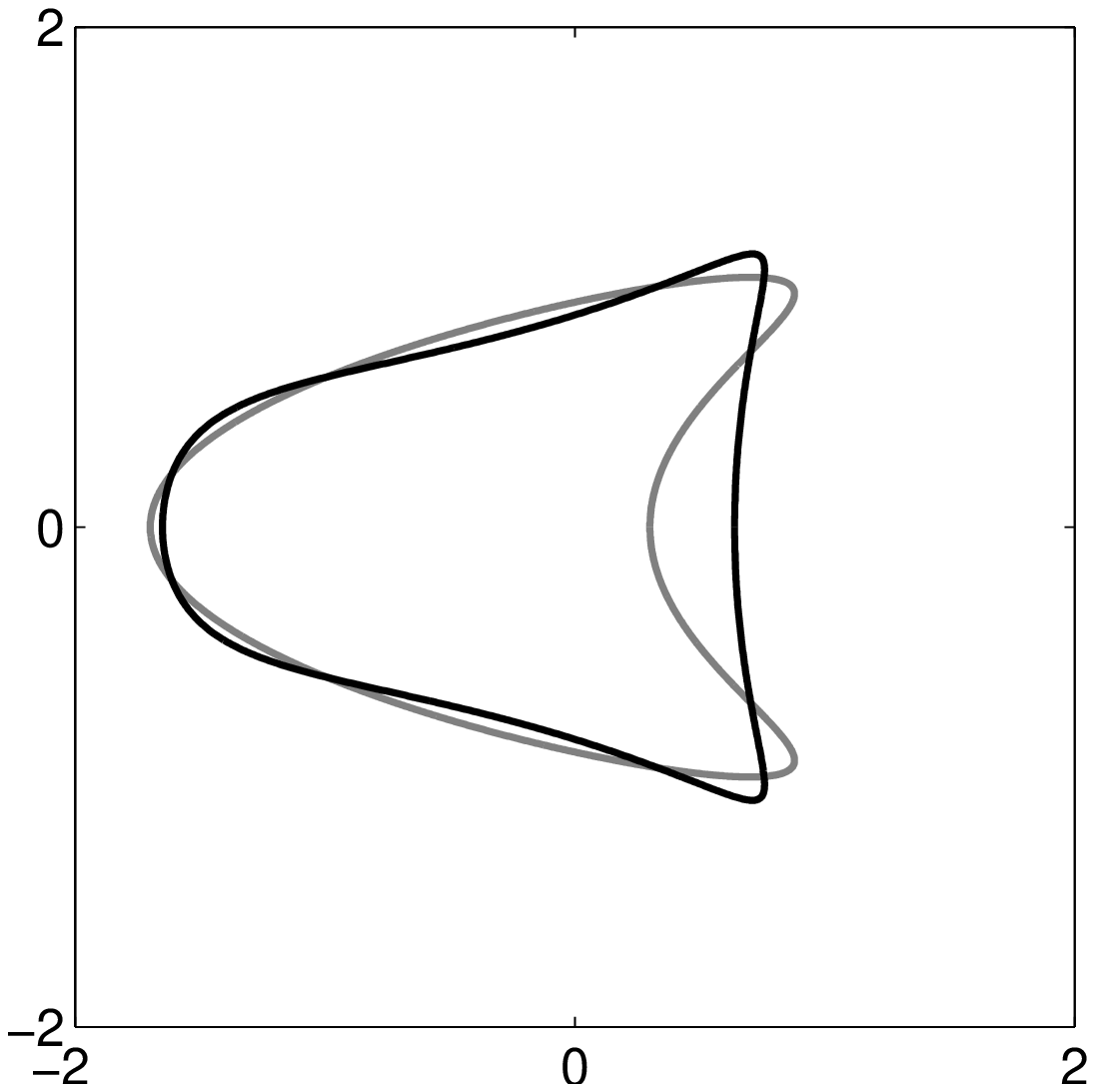,width=3.5cm}\\[.5cm]
\epsfig{figure = 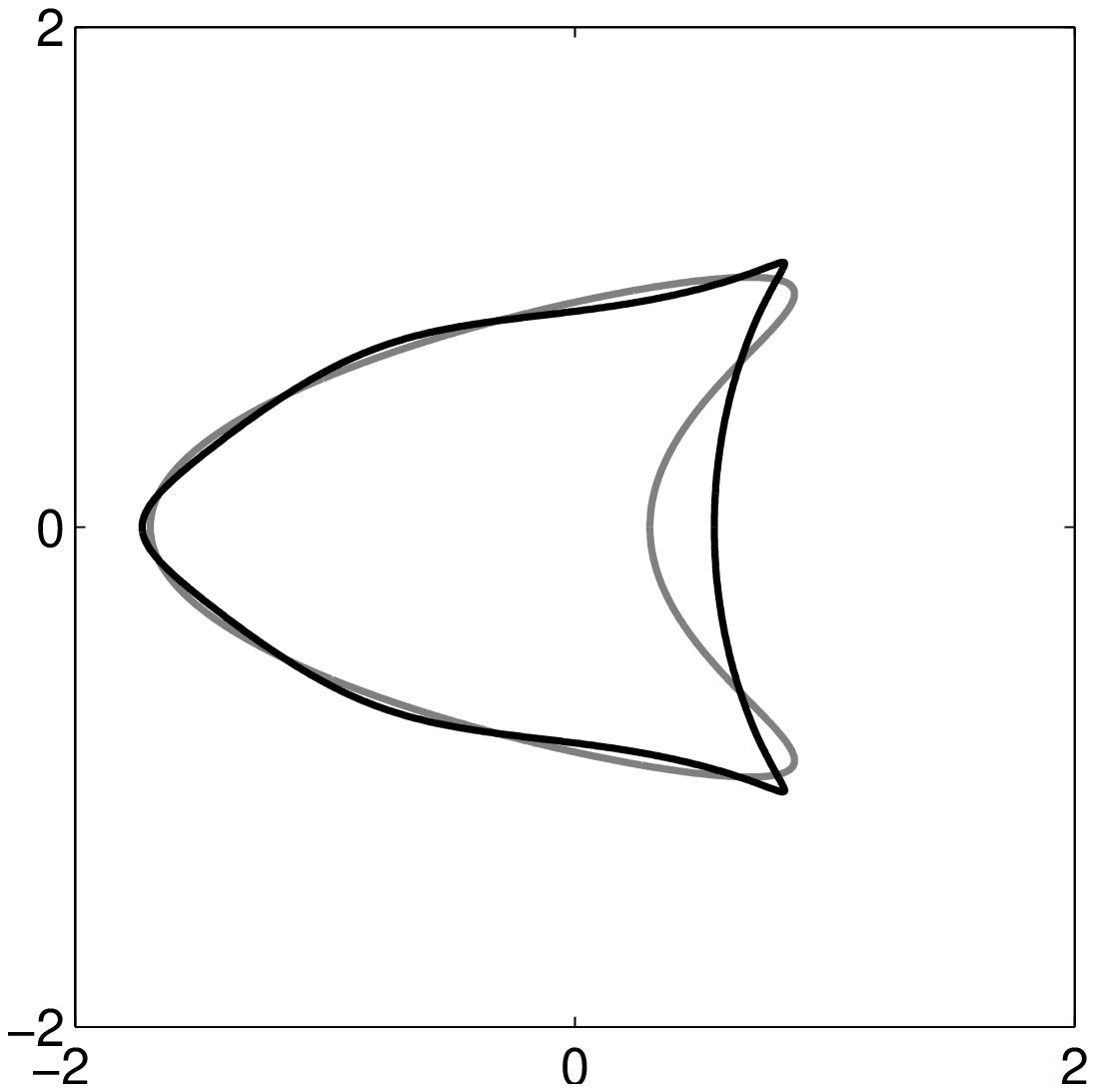,width=3.5cm}\hskip .5cm
\epsfig{figure = 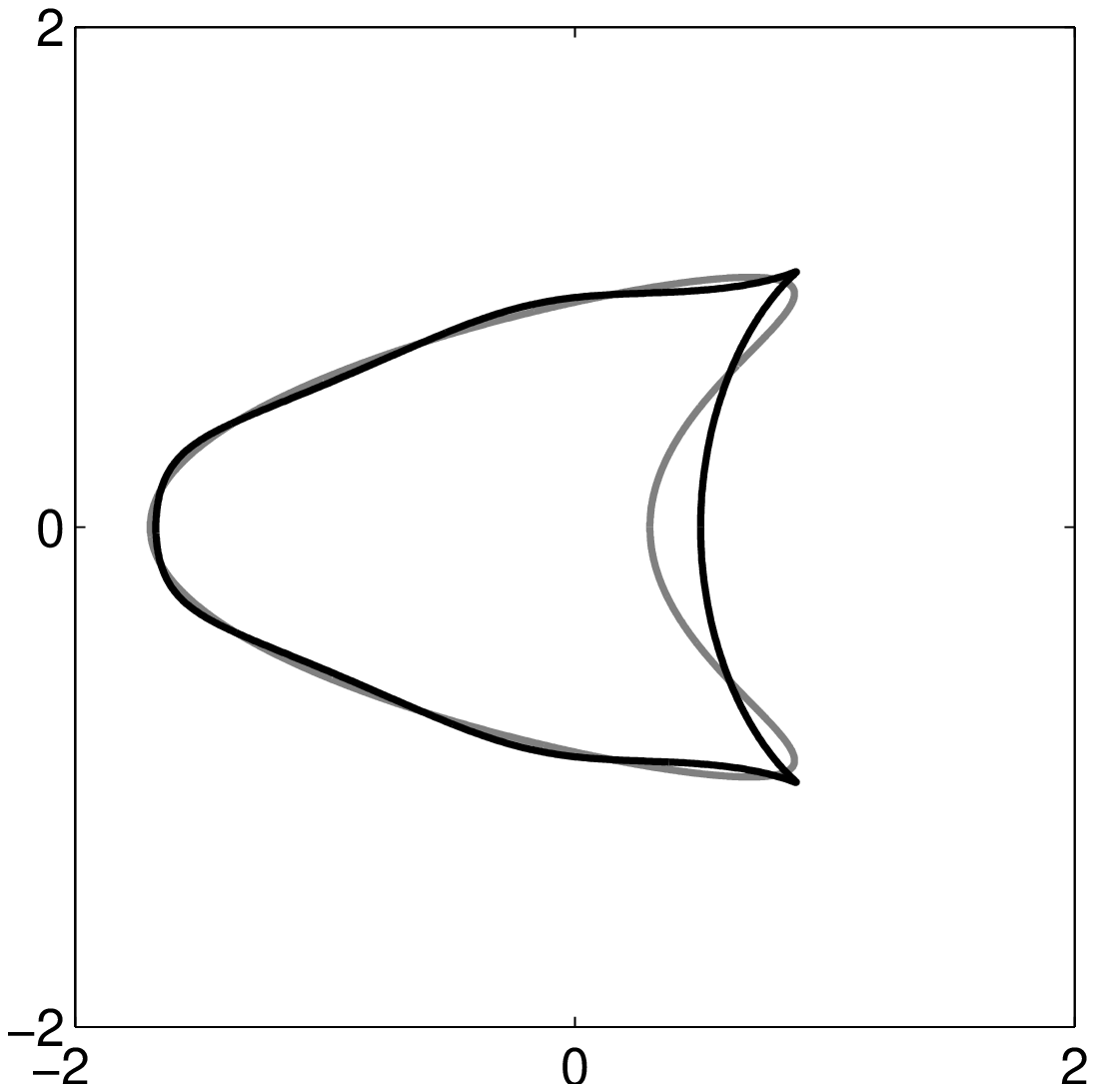,width=3.5cm}\hskip .5cm
\epsfig{figure = 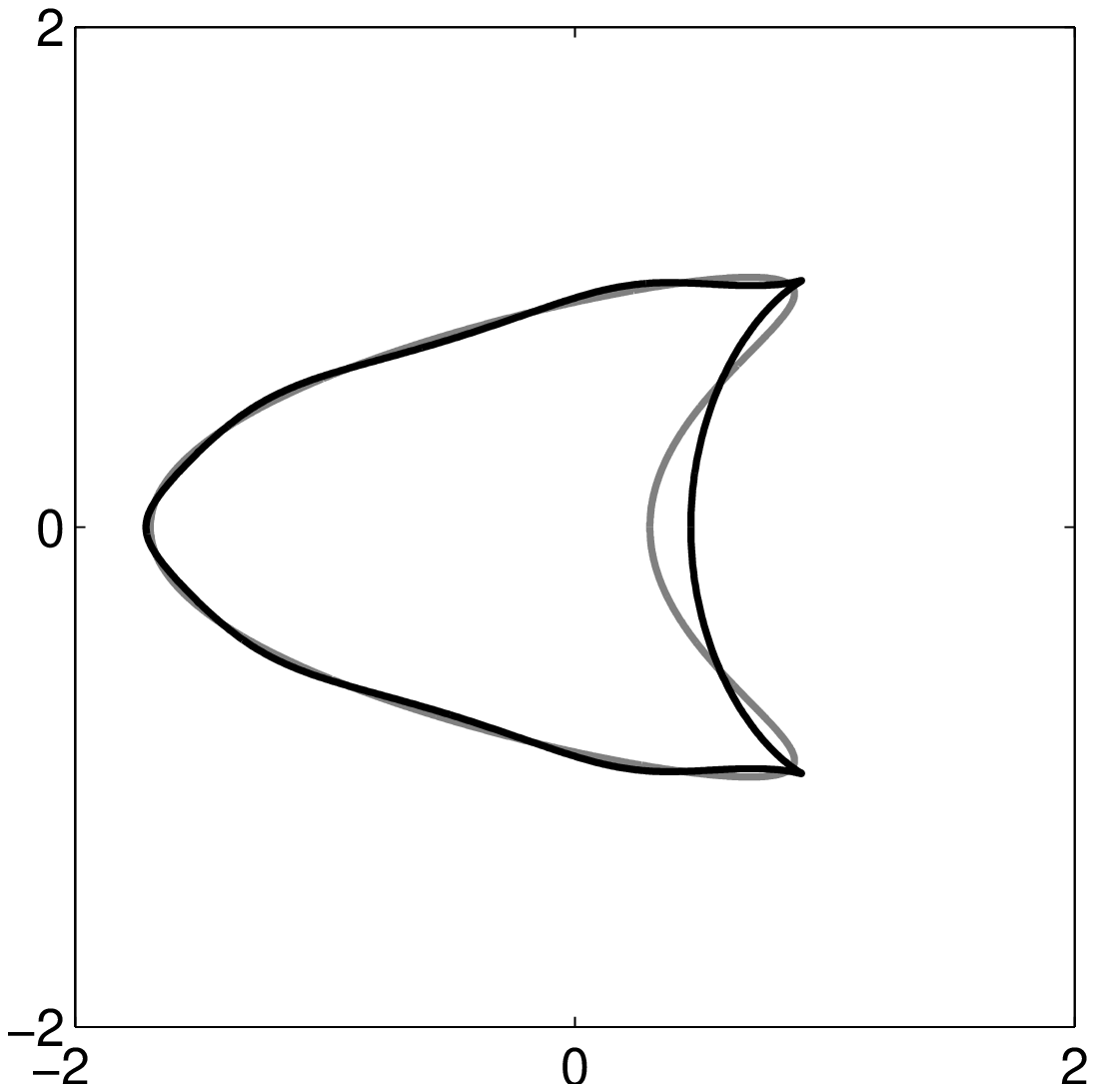,width=3.5cm}
\end{center}
\caption{A kite-shape domain $\Omega$ and $\Phi_N(S^1)$ for $N=1, \ldots, 6$.}\label{kitebig}
\end{figure}

\noindent{\bf Example 2}. Figure \ref{kitebig} shows $\Phi_N(S^1)$ is gradually changing to the boundary of a kite shape domain $\Omega$ as $N$ increases. The computed values of coefficients are presented in Table \ref{table}. The ellipse in the first figure (top left) is called the equivalent ellipse of $\GO$ \cite{book2, bhv}.

\begin{table}
\begin{center}
\begin{tabular}{|c||c|c|c|c|c|c|c|c|c|}
   \hline
$l$&-1&0&1&2&3&4&5&6\\\hline\hline
$\mu_l$& 1.1337 & -0.2415 & 0.1442 &-0.2645&-0.1328&-0.0812&-0.0548& -0.0394\\\hline
\end{tabular}
\end{center}
\caption{The coefficients $\mu_l$, $l\leq 6$, for $\Omega$ in Figure \ref{kitebig}.}\label{table}
\end{table}

\medskip

\begin{figure}[ht!]
\begin{center}
\epsfig{figure = 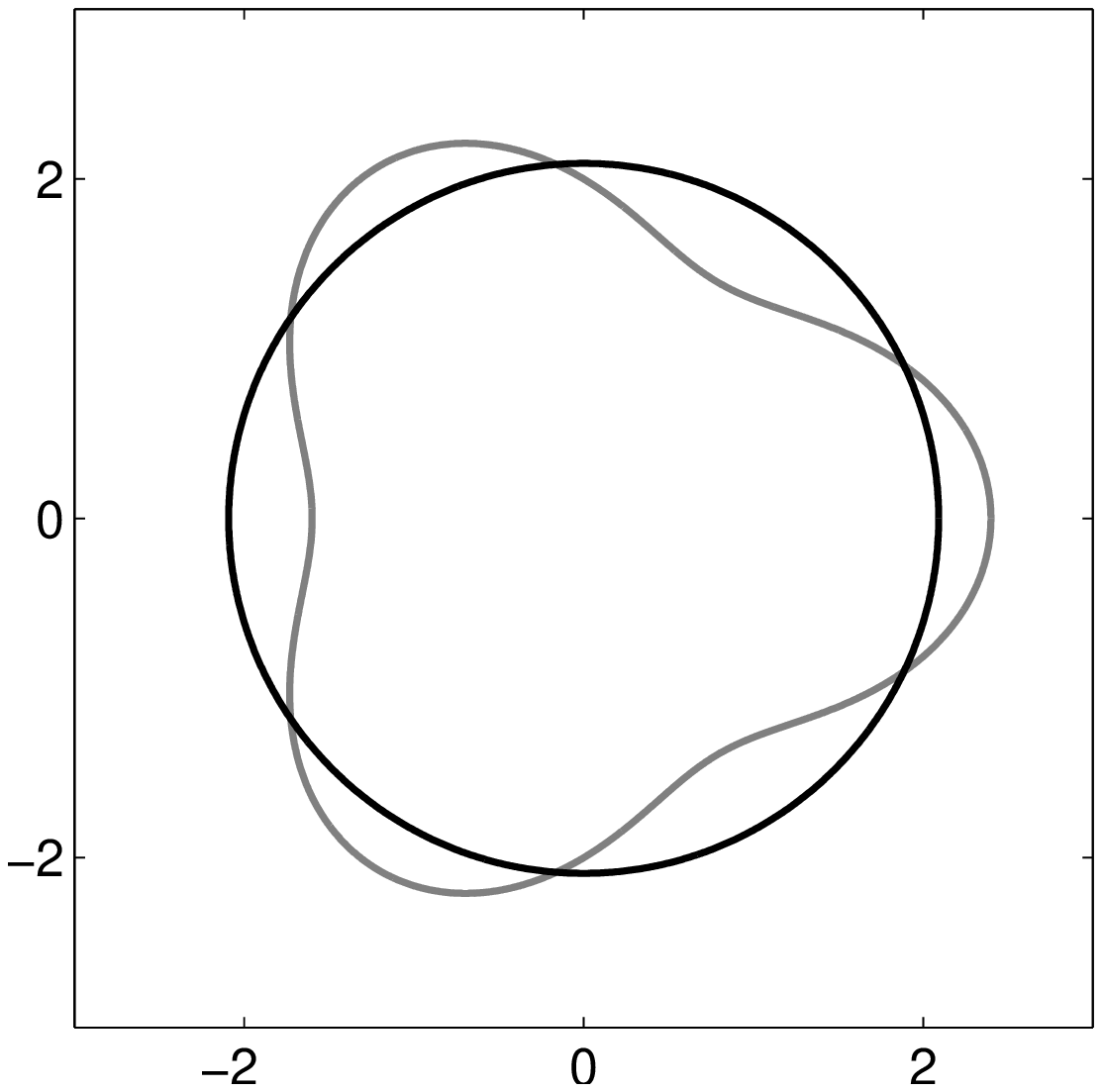,width=3.5cm}\hskip .5cm
\epsfig{figure = 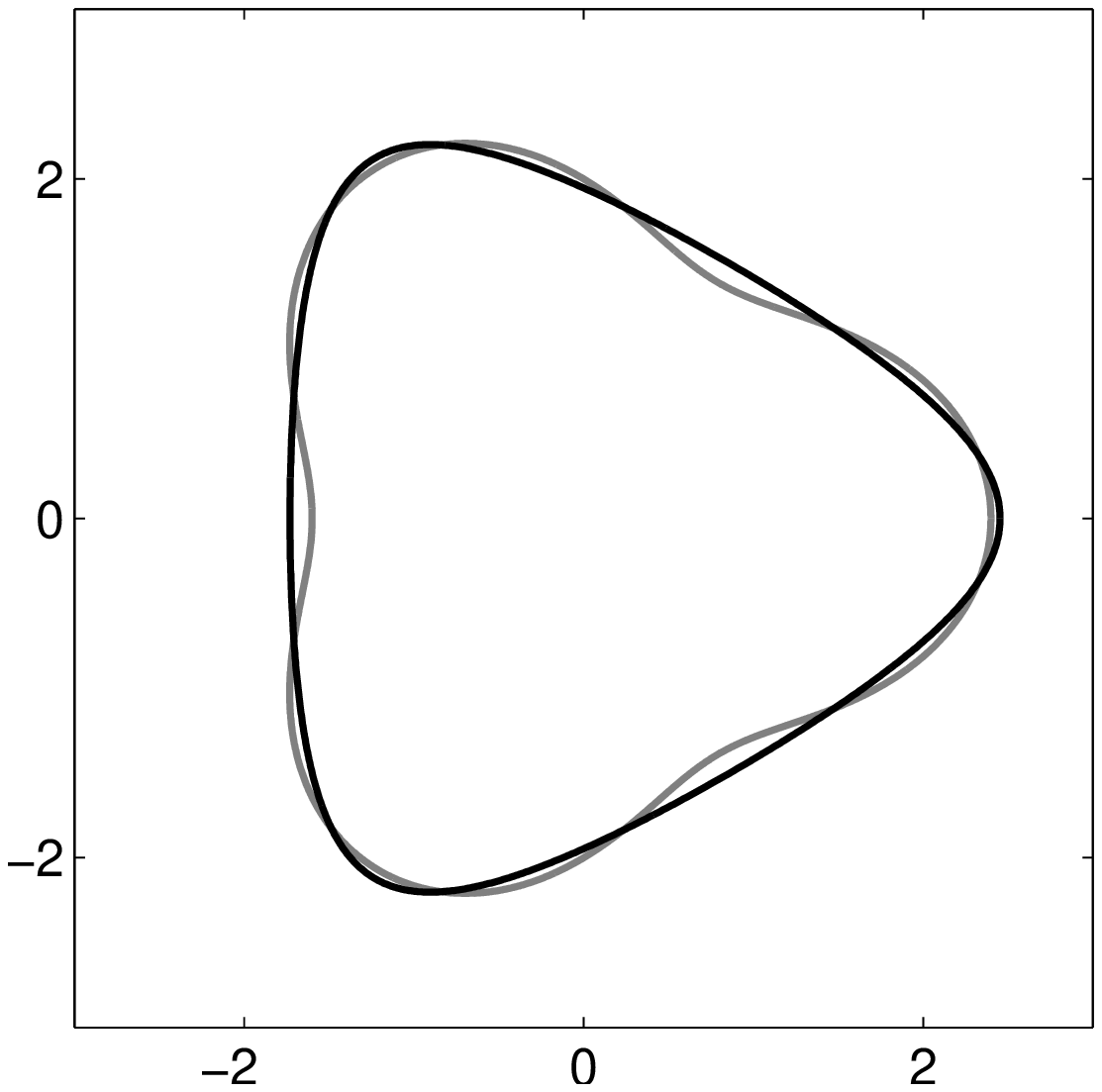,width=3.5cm}\hskip .5cm
\epsfig{figure = 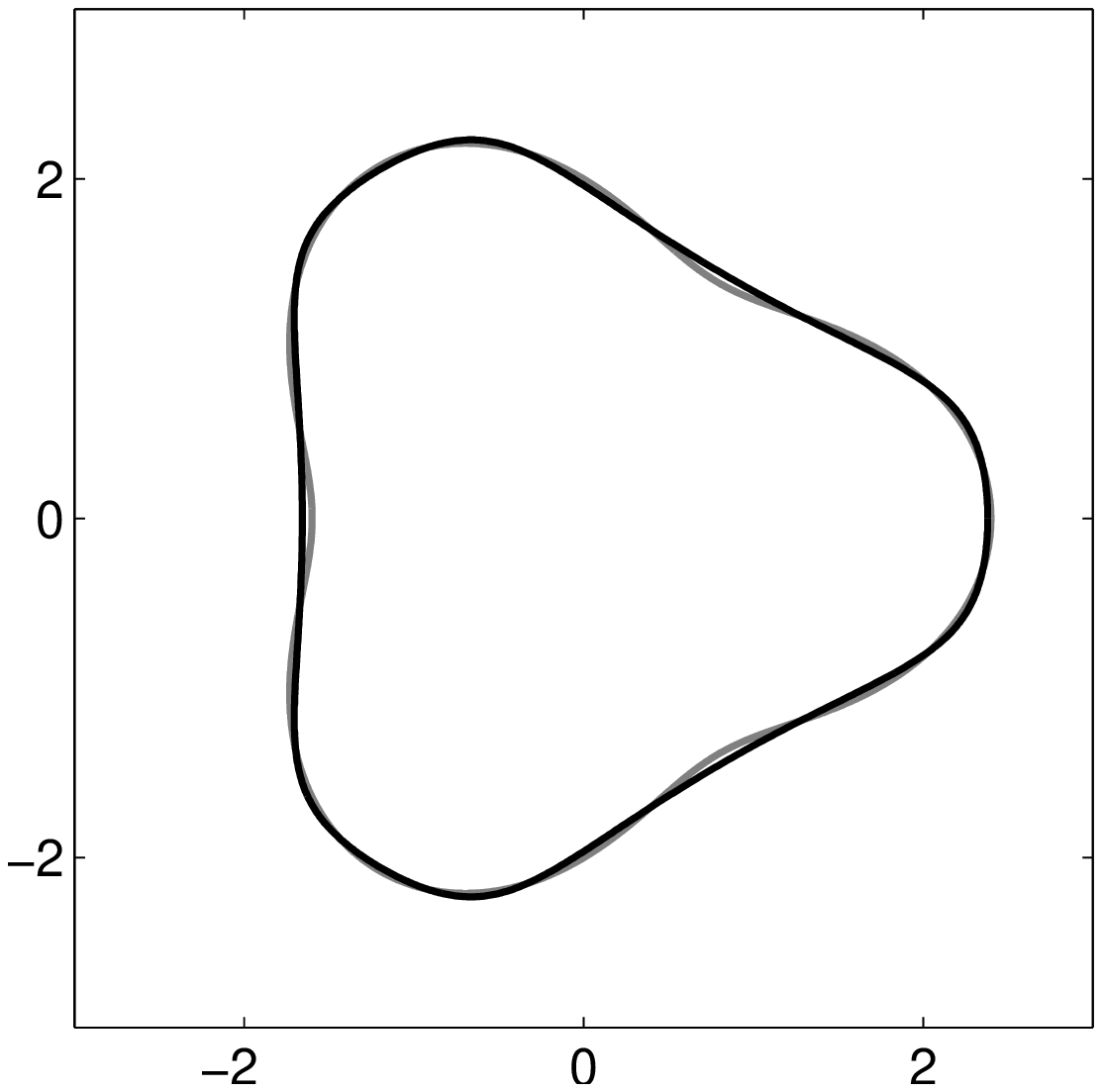,width=3.5cm}\\[.5cm]
\epsfig{figure = 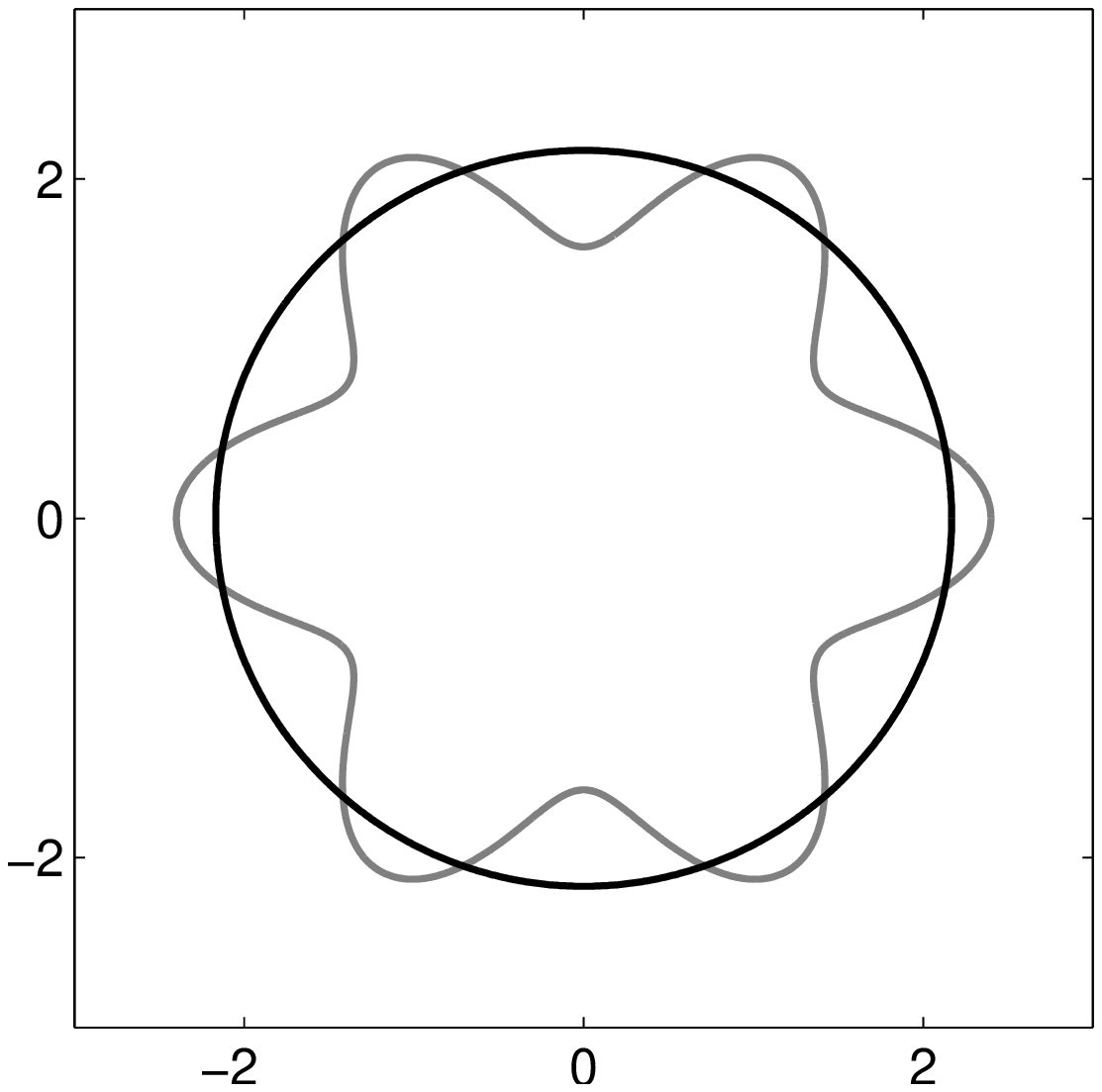,width=3.5cm}\hskip .5cm
\epsfig{figure = 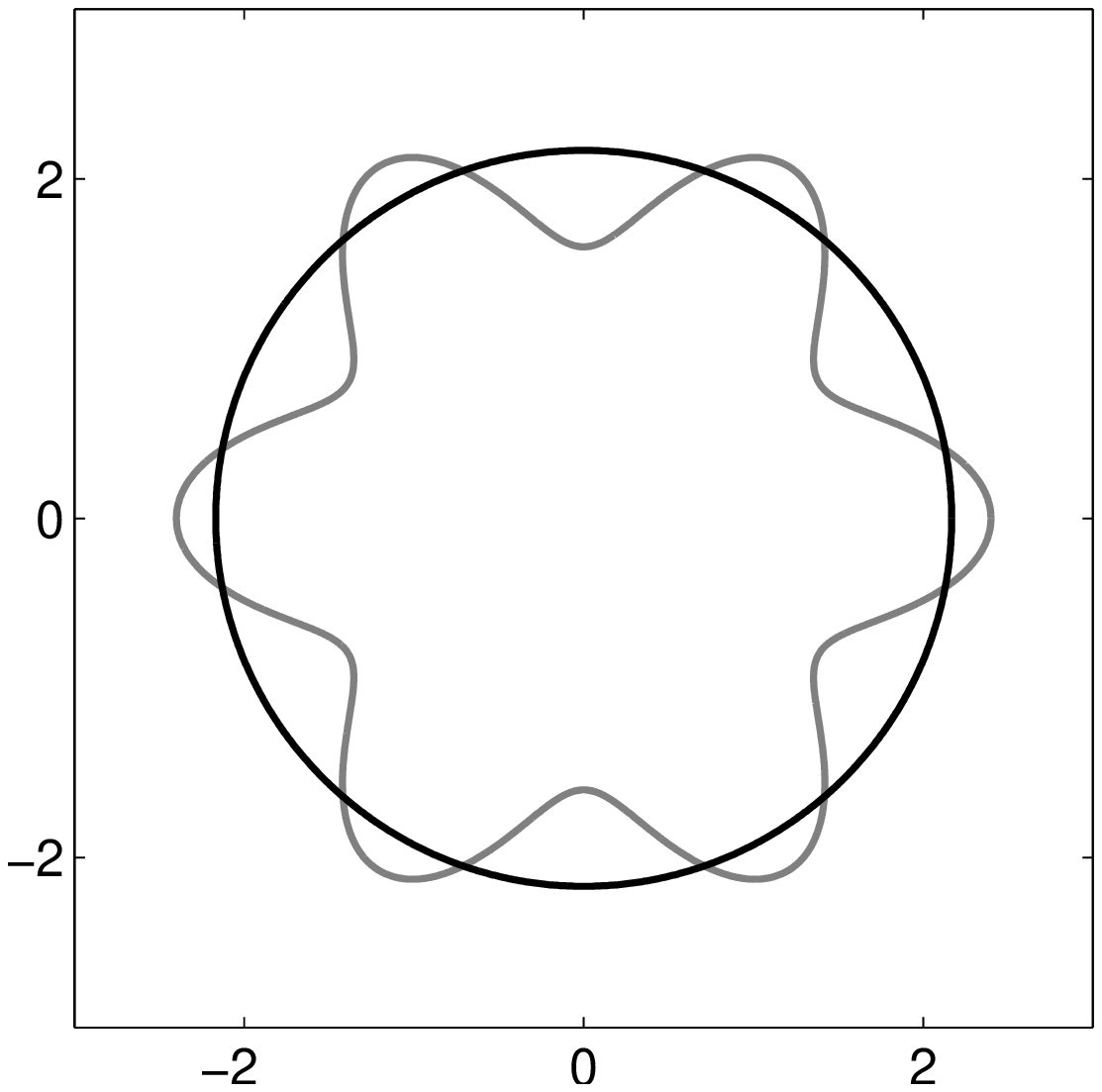,width=3.5cm}\hskip .5cm
\epsfig{figure = 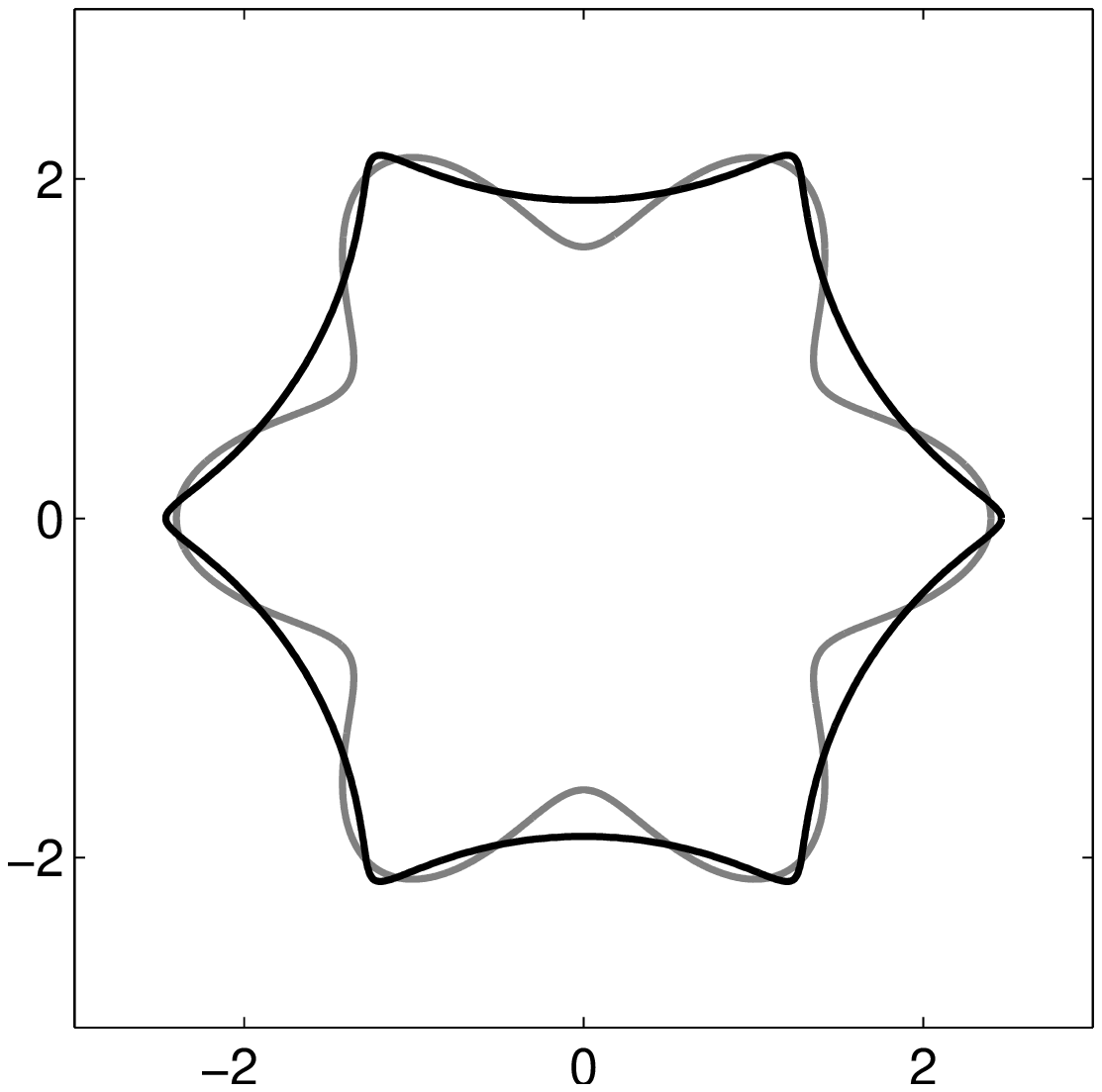,width=3.5cm}
\end{center}
\caption{The gray curve is $\p\Omega$ given by $r=(2+ 0.4\cos(p\theta))$ in polar coordinates for $p=3$ (in the top row) and $p=6$ (in the bottom), and the black is the images of the unit circle under $\Phi_N$ for $N= 1, 2, 5$. }\label{highrow}
\end{figure}

\noindent{\bf Example 3} Figure \ref{highrow} reveals that the boundary with mild oscillation can be recovered by $\Phi_N$ for relatively small $N$, while that with high oscillation requires $\Phi_N$ for higher $N$. This fact was also observed in \cite{aklz}.

%We also provide the examples with various noise levels using  \eqnref{entire01} and \eqnref{constant01} are used as well as \eqnref{c_mu0} and \eqnref{mu_higher}, which is to improve the stability. Adding more equations, the problem becomes overdetermined, and hence the gradient descent method for the suitable energy functional was used.

%%%%%%%%%%%%%%%%%%%%%%%%%%%%%%%%%%%%%%%%%%%%%%%%%%%%%
\section{Further discussion}\label{sec5}
%%%%%%%%%%%%%%%%%%%%%%%%%%%%%%%%%%%%%%%%%%%%%%%%%%%%%

We have derived an explicit connection between GPTs and coefficients of the conformal mapping, and show by numerical examples that first few terms of the conformal mapping approximate the domain quite well.

It is quite interesting to extend results of this paper to construction of conformal mappings of multiply connected domains. We emphasize that GPTs are defined for multiply connected domains as well. In this regard, it is worth emphasizing that only the relations for $n=1$ in \eqnref{entire01} and some partial relations in \eqnref{constant01} are used to derive relation between GPTs and the conformal mapping. So, the relations for $n \ge 2$ and other relations in \eqnref{constant01} provide relations among GPTs. In particular, the equation \eqnref{vanishing_eq} says that  all the terms in $\{\Gg_{m1}^2:~m\ge3\}$ can be calculated by \eqnref{full_set}.  For instance, we obtain
\beq
\Gg_{31}^2=\Gg_{11}^1\Gg_{11}^2 + \frac{(\Gg_{21}^2)^2}{\Gg_{11}^2}.
\eeq
This relation holds only for simply connected domains. For example, if the domain is two disjoint unit disks centered at $(\pm 2, 0)$, then $\Gg_{31}^2=-8.03$ and $\Gg_{11}^1\Gg_{11}^2 + (\Gg_{21}^2)^2/\Gg_{11}^2 = -0.25.$

Note that translation, rotation, and scaling of the domain $\GO$ are expressed as $\Ga \Phi + \Gb$ for some complex numbers $\Ga$ and $\Gb$. So, the quantities $\mu_j/\mu_{-1}$ ($j=1,2,\ldots$) are invariant under translation, rotation, and scaling. In other words, they can be used as shape descriptors in 2D, which can be computed using GPTs. It is worth mentioning that invariant shape descriptors are derived in two and three dimensions using GPTs in \cite{abgj, ackw} and used effectively in a new development of electro-sensing \cite{abgw}.

It is a classical subject to derive optimal bounds for the coefficients of the conformal mapping (see, for instance, \cite{hen} and references therein). In this regards, it is worthwhile to mention the Bieberbach conjecture and its resolution by de Brange \cite{deBrange}. On the other hand, it is an important problem to derive optimal bounds of GPTs. For example, the bounds for the first order GPTs ($\Gg^{1}_{11}$ and $\Gg^{2}_{11}$) are obtained in \cite{cv, lipton}. The relation between GPTs and the conformal mapping obtained in this paper may shed new light on this problem.

\end{document}